\documentclass[11pt]{amsart}
\usepackage{amssymb,latexsym,amsmath}

\newtheorem{prop}{Proposition}

\newcommand{\example}{{\it Example: }}
\newcommand{\remarkr}{{\it Remark: }}
\newcommand{\Problem}{{\it Problem: }}
\newcommand{\claim}{{\it Claim: }}
\newcommand{\Case}{{\it Case: }}

\newcommand{\real}{\mathbb{R}}
\newcommand{\diag}{\text{Diag}}
\newcommand{\Trace}{\text{Trace}}
\newcommand{\Rank}{\text{Rank}}
\newcommand{\Orbit}{\text{Orbit}}
\newcommand{\Tan}{\text{Tan}}

\newcommand{\birkhoff}{Birkhoff and MacLane(1953) }
\newcommand{\driessel}{Driessel(2004)}

\begin{document}

\title[Low Rank Approximation]
{A Globally Convergent Flow for \\
Computing the Best Low Rank Approximation \\
of a Matrix  }
\author{Kenneth R. Driessel}
\address{Mathematics Department\\
Iowa State University}
\date{\today}         

\keywords{equivalence, singular values,
gradient, flow, quasi-projection, 
Eckart-Young theorem,
approximation, Riemannian manifold}

\subjclass[2000]
{Primary: 15A03 Vector spaces,linear dependence, 
rank;
Secondary: 15A18 Eigenvalues, singular values,
and eigenvectors}
                                
\maketitle

\newpage

\section*{Abstract}

We work in the space of $m$-by-$n$
real matrices with the Frobenius inner
product. Consider the following problem:
\smallskip

\Problem: Given an m-by-n real matrix
A and a positive integer k, find the m-by-n 
matrix with rank k that is closest to A. 
\smallskip

I discuss a rank-preserving  differential 
equation (d.e.)  
which solves this problem. If X(t) is a solution
of this d.e., then the distance between $X(t)$ 
and $A$ decreases as t increases; 
this distance function 
is a Lyapunov function for the d.e. 
If $A$ has distinct 
positive singular values (which is 
a generic condition) then this d.e.
has only one stable equilibrium point.
The other equilibrium points are finite
in number and unstable.
In other words, the basin of attraction
of the stable equilibrium point on the 
manifold of matrices with rank $k$ 
consists of almost all matrices.
This special equilibrium point is 
the solution of the given problem.
Usually constrained optimization problems
have many local minimums (most of which 
are undesirable). So the constrained
optimization problem considered here
is very special. 

\section*{Table of Contents}

\begin{itemize}

\item
Introduction
\smallskip

\item
Setting Up the Differential Equation
\smallskip

\item
Properties of the Differential Equation
\smallskip

\item
Acknowledgements

\item
Appendix: The Frobenius Inner Product
\smallskip

\item
References

\end{itemize}

\newpage
\section{Introduction}

We work in the space $\real^{m\times n}$
of $m$ by $n$ real matrices with
the ``Frobenius'' (or ``euclidean'')
inner product. 
(In an appendix we
review the definition and
elementary properties of this 
inner product.) We consider the 
following problem:

\smallskip
\Problem {\bf Low rank approximation.}
Given a matrix $A$ in 
$\real^{m\times n}$ and a positive
integer $k$, find the matrix
with rank $k$ which is closest to $A$. 
\smallskip

This problem is closely connected
with the singular value decomposition
of matrices. If $A=UDV^T$ where 
$U$ is an $m\times m$
orthogonal matrix, $V$ is an $n\times n$
orthogonal matrix, and
$D=\diag(\sigma_1\geq\sigma_2\geq\dots \geq\sigma_n)$
is a diagonal matrix,
then the product $UDV^T$ is called the
{\bf singular value decomposition of $A$}.
The diagonal entries of $D$ are the
{\bf singular values of $A$}. It is
well-known that every matrix has a
singular value decomposition. (See, for
example, Horn and Johnson(1985) or Demmel(1997).)

The following
result provides a solution of the 
low rank approximation problem:

\begin{prop}
Assume $m\geq n$ and that the rank of $A$
is greater than the positive integer
$k$. Let $A=UDV^T$
be the singular value decomposition
of $A$. Let 
$D^\prime:=\diag(\sigma_1,\sigma_2,\dots,\sigma_k,
0,\dots,0)$ and let
$A^\prime:=UD^\prime V^T$. Then $A^\prime$
is the matrix with rank $k$ which most
closely approximates $A$ in the
Frobenius norm. 
\end{prop}

This result is often called the Eckart-Young
theorem. The result appeared in
Eckart-Young(1936). However, Stewart(1993)
points out that it was known earlier.

For a textbook proof of this proposition,
see, for example, Horn and Johnson(1985) Section 7.4
``Examples and applications of the singular
value decomposition". I 
present an alternative proof of this 
result in this paper. In particular, I 
discuss a quasi-gradient differential equation
which computes the solution of the given
problem.  This proof provides more
information than other proofs. In 
particular, it shows that (generically)
$A^\prime$ is the unique local minimum for the
low rank approximation problem. In other
words, the basin of attraction of this
matrix consists of almost all matrices 
on the surface of matrices with rank $k$.  

The proposition suggests a way to compute
the solution of the low rank
approximation problem: Compute the
singular value decomposition of $A$
then compute the approximation $A^\prime$.
This procedure is obviously inefficient:
Why compute all the singular values of $A$
if we only need the largest ones for 
the solution? We shall see that the
differential equation is more 
``economical'' since its flow is on
the manifold $\Rank(k)$ of matrices
with rank $k$. If $k$ is small then
this manifold has dimension much
smaller than the dimension of 
$\real^{m\times n}$. I hope that
this differential equation can be 
used to design an efficient algorithm
for low rank approximation.

Since the 1980's there has been significant
work with flows on manifolds of matrices. 
In particular, 
during the 1980's, there was considerable
interest in continuous analogues of 
the QR algorithm for computing 
eigenvalues of matrices. 
The connection between the QR algorithm
and the Toda flow was discovered by
Symes about 1980. 
For more on this connection,  
see, for example,
Symes(1980a,1980b,1982), 
Deift, Nanda and Tomei (1983),
Nanda (1982,1985), Chu(1984), and 
Watkins (1984a,1984b).
There are now also textbook descriptions of 
this connection: See, for example, Demmel (1997).
For some other
flows on matrices, see 
Chu(1986a,1986b), Chu and Driessel (1990),
Helmke and Moore(1995), 
Driessel(2004), Driessel and Gerisch(2007)
and the works cited
in these references.  

Flows on manifolds of matrices 
are interesting not just
because of their connections with 
computation, but also for the insight they
provide into the geometry of the
manifolds of interest. 
This is the main idea in Morse theory.
Let me say more about such geometric
insights. 
Let $W$ be a real vector space with an
inner product $\langle \cdot,\cdot\rangle$. Let
$S$ be a subset of $W$ and let 
$f:S\to \real$ be a real-valued
function on $S$. Then $m\in S$ is
a {\bf local minumum} of $f$ on $S$,
if there is a neighborhood $N$ of $m$
such that $f(m)$ is a minimum of
$f$ on $N$. I say that $f:S\to\real$ has
the {\bf  unique local minimum property}
if $f$ is bounded below and has a unique
local minimum. In this case the local 
minimum is also the global minimum.
Usually an optimization
problem has numerous (mostly undesirable)
local minimums. An optimization 
problem with the unique local minimum
property is an especially nice 
optimization problem. 

Here are a few examples. 
Let $S$ be a convex set in $\real^n$
with the euclidean inner product;
let $a$ be a point in $\real^n$; let
$f_a:S\to \real$ be defined by
to be the square of the distance
from $s\in S$ to $a$: 
$f_a(s):=\langle a-x,a-x \rangle$;
this function has the unique local
minimum property for all $a$. 
Let  $S$ be a circle in the euclidean
plane $\real^2$ and, for a point
$a$ in $\real^2$ let $f_a:S\to \real$
be the square of the distance from 
$s$ to $a$; this function has the
unique local minimum property
unless $a$ is the center of the circle. 

Here is another example. Consider the
following problem:

\Problem {\bf Approximation with 
spectral constraint.} Given an
$n\times n$ symmetrix matrix and
real eigenvalues 
$\lambda_1,\lambda_2,\dots,\lambda_n$,
find the matrix with these eigenvalues
that is closest (in the Frobenius
norm) to $A$.

Chu and Driessel(1990) studied this
problem. They showed (by means
of a ``gradient'' flow) that
it satisfies the unique local minimum
property if the eigenvalues
$\lambda_i$  are distinct.

Let $\Rank(k)$ denote the set of 
matrices in $\real^{m\times n}$
with rank $k$. Let $A$ be a matrix
in $\real^{m\times n}$. Define
the function $f_A:\Rank(k)\to \real$
by $f_A(X):=(1/2)\langle A-X,A-X \rangle$
where $\langle \cdot,\cdot \rangle$
is the Frobenius inner product. 
In this paper I show  
that if $A$ has distinct positive
singular values then the function 
$f_A$ has a finite number of 
critical values only one of which
is a local minimum and hence $f_A$ has
the unique local minimum property.

\remarkr Helmke and Shayman(1985),
in Theorem 4.2(ii), say that $f_A$
has a finite number of critical points
if and only if $m=n$ and $A$ has $m$
distinct nonzero singular values.
The results that I present here
show that the condition $m=n$ is
not necessary.
 
 {\it Contents summary:}
 In the section with the title
 ``Setting up the differential equation'',
 I review the differential geometry 
 associated with the rank approximation
 problem. (This material appears in
 Helmke and Moore(1995) and in 
 Helmke and Shayman(1995). I include it
 to make this paper more self-contained.)
 I also describe the quasi-projection
 operator associated with this problem. 
 (For more on such operators see 
 Driessel(2004).) 
 In the section with title ``Properties
 of the differential equation'', I show 
 that this differential equation has the
 convergence properties asserted above. 
 (This differential equation appears in
 Helmke and Moore(1995) and 
 Helmke and Shayman(1995) but they
 derive it in a more complicated way
 than I do.
 Their discussion of its equilibrium
 points is not very clear.  
 They do not classify the equilibrium
 points. They do not discuss basins of 
 attraction.)
 
 The only prerequisites for understanding
 (almost all) of this paper are a basic
 knowledge of differential equations
 (see, for example, Hirsch and Smale(1974))
 and basic differential geometry
 (see, for example, Thorpe(1979)).
 
\newpage
\section{Setting up the differential 
equation}

In this section we view the sets 
$\Rank(k)$ of matrices 
with fixed ranks $k$ as parameterized 
surfaces in the space $\real^{m\times n}$.
We compute the tangent spaces of these constant
rank surfaces. Then 
we define a ``quasi-projection'' map which
can be used to transform vector fields
in $\real^{m\times n}$ into vector
fields tangent to these surfaces.
We also define an objective
function associated with the 
constrained optimization problem of
interest and we compute its gradient. 
Finally, we use the quasi-projection 
map to convert this gradient
vector field into one which is tangent
to the constant rank surfaces. 

Let $Gl(m)$ denote the general linear
group of $m$ by $m$, invertible, real
matrices. Recall (see, for example, 
\birkhoff) that two matrices $X$ and
$Y$ in $\real^{m\times n}$ are
{\bf equivalent} if there exist
matrices $G\in Gl(m)$ and $H\in Gl(n)$
such that $Y=GXH^{-1}$. Also recall
that every matrix $M$ in 
$\real^{m\times n}$ is equivalent to 
a diagonal matrix $D$ with ones and
zeros on its main diagonal. The number
of ones equals the rank of $M$.

We can use the groups $Gl(m)$ and 
$Gl(n)$ to ``parameterize'' the 
matrices with rank $k$ as follows.
We use the following group action:
$$
Gl(m)\times Gl(n)\times \real^{m\times n}
\to \real^{m\times n}:
(G,H,X)\mapsto GXH^{-1}.
$$
For $M\in\real^{m\times n}$, we use $\Orbit(M)$
to denote the orbit of $M$ under this group
action; in symbols, 
$$
\Orbit(M):=\{GMH^{-1}:G\in Gl(m),
H\in Gl(m)\}.
$$
Let 
$$
\Rank(k):=\{X\in\real^{m\times n}:
\Rank(X)=k\}.
$$
The following proposition summarizes the
comments given above.

\begin{prop}
Let $k$ be a positive integer and let $K$
be any matrix with rank $k$. Then
the set of matrices with rank $k$
is the same as the orbit of $K$
under the given group action; in symbols,
$$
\Rank(k)=\Orbit(K).
$$
\end{prop}

For $B\in\Orbit(K)$, I use 
$\Tan.\Orbit(K).B$ to denote the
space tangent to $\Orbit(K)$ at $B$. 

\begin{prop}
Let $K$ and $B$ be matrices in 
$\real^{m\times n}$ with $B$ on the orbit
of $K$. 
Then the space tangent to the orbit of $K$
at $B$ is given by
$$
\Tan.\Orbit(K).B =
\{XB+BY: X\in\real^{m\times m},
Y\in\real^{n\times n}\}.
$$
The dimension of this tangent space
is $k^2+k(m-k)+k(n-k)$
where $k:=\Rank(K)$.
\end{prop}

\begin{proof}
We simply compute the derivative of 
the parameterizing map. We have
(where $D$ denotes the derivative
operator)
\begin{align*}
&D((G,H)\mapsto GBH^{-1}).(I,I).(X,Y) \\
&= (D(G\mapsto GB).I.X)\cdot(I^{-1})+
(IB)\cdot(D(H\mapsto H^{-1}).I.Y) \\
&= XB + BY
\end{align*}
since
$$
D(H\mapsto H^{-1}).C.W = - C^{-1}W C^{-1}.
$$
(I sometimes use dots for function
evaluation in order to reduce the 
number of parentheses. I also use
association to the left.)

Since the orbit is a homogeneous space, it
looks the same at all its points. Consequently,
we can compute the dimension of the
tangent space at any convenient point.
For example, we can do the computation
at the diagonal matrix 
with exactly $k$ ones on its diagonal
and zeros elsewhere; in particular,
we can take  
$B:=\diag(1^{\times k},0^{\times (l-k)})$
where $l$ is the minimum of $m$ and $n$. 
\end{proof}

For $B$ on the orbit of $K$, we consider the
following linear map:
$$
L_B:=
\real^{m\times m}\times\real^{n\times n}
\to \real^{m\times n}:
(X,Y)\mapsto XB+BY.
$$
Note that the range of this map 
equals the space tangent to the
orbit of $K$ at $B$. We compute the
adjoint $L_B^*$ of this map.

\begin{prop} {\bf Adjoint of the 
tangent space map.}
Let $B$ be on the orbit of $K$. Then
the adjoint $L_B^*$ of the linear map $L_B$
is the map
$$
\real^{m\times n}\to 
\real^{m\times m}\times\real^{n\times n}:
Z\mapsto (ZB^T,B^TZ).
$$
\end{prop}

\begin{proof}
We have
\begin{align*}
\langle L_B(X,Y),Z\rangle 
&= \langle XB+BY,Z \rangle
= \langle X,ZB^T \rangle + \langle Y,B^TZ\rangle \\
&= \langle (X,Y),(ZB^T,B^TZ) \rangle.
\end{align*}
Here we have used the ``product'' inner
product on the space 
$\real^{m\times m}\times\real^{n\times n}$ which
is defined in terms of the Frobenius
inner product by:
$$
\langle (X_1,Y_1),(X_2,Y_2) \rangle :=
\langle X_1,X_2 \rangle +
\langle Y_1,Y_2 \rangle .
$$
\end{proof}

I call the composition $L_B\circ L_B^*$
a ``quasi-projection'' map. 
We can use this operator to transform
vector fields on $\real^{m\times n}$
into ones which are tangent to the
orbits of interest.
For more on the use of quasi-projections
 see \driessel.

For $A$ in $\real^{m\times n}$, we define 
the {\bf objective function $f:=f_A$ 
determined by} $A$ as the following function:
$$
\real^{m\times n}\to \real:
X\mapsto (1/2)\langle X-A,X-A \rangle.
$$
In other words, $f_A(X)$ is one half
the square of the distance from $X$ to $A$.
\begin{prop}{\bf Gradient of the objective
function.}
Let $A$ and $B$ be matrices in 
$\real^{m\times n}$. The gradient of 
the objective function $f_A$ at $B$
is $B-A$; in symbols, 
$$
\nabla f_A(B) = B-A.
$$
\end{prop}

\begin{proof}
We simply compute the derivative of $f:=f_A$:
For $X$ in $\real^{m\times n}$, we have
\begin{align}
Df.B.X 
&= D(X\mapsto (1/2)\langle X-A,X-A\rangle).B.X \\
&= \langle B-A, D(X\mapsto X-A).B.X \rangle
= \langle B-A,X \rangle.
\end{align}
\end{proof}
 
 We have a gradient vector field on
 $\real^{m\times n}$ defined by
 $X\mapsto \nabla f_A (X)=X-A$. 
 But this vector field is generally
 not tangent to the constant rank
 surfaces. In other words, the corresponding
 differential equation $X^\prime=\nabla f_A(X)$
 does not preserve rank. We want to adjust the
 gradient vector field so that the 
 corresponding vector field does 
 preserve rank. We can use the
 quasi-projection map to do so. 
 
 We now compute the quasi-projection
 of the negative gradient onto the 
 tangent space. For $B$ on the orbit of 
 $K$, we have
 \begin{align*}
 (L_B \circ L_B^*) (-\nabla f_A (B))
 &= L_B((A-B)B^T,B^T(A-B)) \\
 &=(A-B)B^TB+BB^T(A-B).
 \end{align*}
In the next section we shall use this 
formula to define a vector field
on the space $\real^{m\times n}$.
We shall then see that the 
corresponding differential equation
provides a solution of the constrained
optimization problem of interest. 

\newpage
\section{Properties of the Differential
Equation}

In the last section we saw how to  adjust
the gradient vector field determined by
the objective function so that the
resulting vector field is tangent to
constant rank submanifolds. We now 
use that quasi-gradient vector field 
to define a differential equation.

Using the results of the last section
we define the vector field 
$F$ on $\real^{m\times n}$
as follows:
$$
F(X):= (L_X\circ L_X^*)(A-X)
= (A-X)X^T X + X X^T (A-X).
$$
We consider the differential equation
associated with this vector field:
\begin{equation*}
X^\prime = F(X). \tag{*}
\end{equation*}

We shall see later that the 
solutions of this differential 
equation are defined for all time. 
In particular we shall see that 
the solutions do not blow up. We shall
also see that they converge.

Note that this differential 
equation is clearly rank preserving
since the vector $F(X)$ is tangent 
to the space $\Rank(X)$ at $X$. 
The following proposition provides
a more concrete argument. 
(In the following analysis, 
we shall only
use the fact that the differential 
equation preserves rank. We 
shall not use the other
assertions of this result.)

\begin{prop}{\bf Rank preserving.}
Let $X(t)$ be the solution
of the initial value problem
\begin{equation*}
X^\prime = F(X),\quad X(0)=K. 
\end{equation*}
Let $G(t)$ and $H(t)$ be 
solutions of the following 
initial value problems:
\begin{align*}
G^\prime &= (A-X)X^TG
,\quad G(0)=I \\
H^\prime &= -X^T(A-X)H,\quad H(0)=I.
\end{align*}
Then $X(t)=G(t)KH(t)^{-1}$
and the rank is invariant. 
\end{prop}

\remarkr The differential equation
for $G$ is determined by a tangent 
vector field on $Gl(m)$ and the
differential equation for $H$ is determined
by a tangent vector field on $Gl(n)$.
Note that the expressions $(A-X)X^T$ 
and $X^T(A-X)$ appear in the expression
defining the vector field $F(X)$. 

\begin{proof}Let $Z(t):= G(t)^{-1}X(t)H(t)$. 
Note $Z(0)=K$ and
\begin{align*}
Z^\prime 
= &-G^{-1}G^\prime G^{-1}XH
+G^{-1}X^\prime H + G^{-1}X H^\prime \\
= &-G^{-1}(A-X)X^TGG^{-1}XH
+ G^{-1}((A-X)X^TX+XX^T(A-X))H \\
&-G^{-1}XX^T(A-X)H \\
=\ &0.
\end{align*}
Hence $Z(t)=K$ for all $t$.
\end{proof}
 
 The following proposition says that 
 for any solution $X(t)$ of the differential 
 equation (*), the distance between $X(t)$
 and $A$ decreases.
 
 \begin{prop}{\bf Lyapunov function.}
 The objective function $f_A$ is a
 Lyapunov function for the differential
 equation (*).
 \end{prop}
 
 \begin{proof}
 Let $X(t)$ be any solution of (*). 
 To simplify the notation, let
 $f:=f_A$ and $L:=L_X$. We have
 \begin{align*}
 (d/dt)(f(X(t))
 &= (1/2)(d/dt) \langle X-A,X-A \rangle \\
 &= \langle X-A, X^\prime \rangle \\
 &= \langle X-A, -(L\circ L^*)(X-A) \rangle \\
 &= - \langle L^*(X-A), L^*(X-A) \rangle \leq 0.
 \end{align*}
 \end{proof}
 
 \begin{prop}
 The solutions of the differential equation (*)
 are defined for all positive times. 
 \end{prop}
 
 \begin{proof}
 Let $X(t)$ be a solution of the 
 differential equation. By the 
 last proposition the distance between
 $A$ and $X(t)$ decreases as $t$ 
 increases. Hence the solution 
 remains in the closed ball
 with radius $\Vert A-X(0) \Vert$
 centered at $A$. Since this ball
 is compact the solution cannot
 blow up.
 \end{proof}
 
 \begin{prop} {\bf Equilibrium conditions.}
 Let $E$ be an element of $\real^{m\times n}$.
 Then the following conditions are equivalent:
 \begin{itemize}
 \item
 [(i)] $E$ is an equilibrium point of 
 the differential equation (*).
 \item
 [(ii)] $E$ satisfies the equations
 $$
 AE^T=EE^T,\quad E^TA=E^TE.
 $$
 \item
 [(iii)] $A-E$ is orthogonal to the space
 tangent the orbit of $E$ at $E$.
 \item
 [(iv)] $E$ is a critical point of 
 the objective function $f_A$.
 \end{itemize}
 \end{prop}
 
 \begin{proof}(i) implies (ii): Let $E$ be an
 equilibrium point of (*). Then (by the
 proof of the Lyapunov proposition) 
 $$
 0 = L_E^*(A-E) = ((A-E)E^T,E^T(A-E)).
 $$
Hence $(A-E)E^T=0$ and $E^T(A-E)=0$.

(ii) implies (i): Assume the $E$ satisfies
the given equations. Then we have
$L_E^*(A-E)=0$
and hence $(L_E\circ L_E^*)(A-E)=0$. 

(ii) implies (iii): Assume that $E$ satisfies
the given equations. Then for any $X$ in 
$\real^{m\times m}$
and $Y$ in  $\real^{n\times n}$, we have
$$
\langle A-E, XE+EY \rangle =
\langle (A-E)E^T, X\rangle +
\langle E^T(A-E), Y \rangle =0.
$$
 
 (iii) implies (ii): Assume that $A-E$
 is orthogonal to the tangent space. Then,
 for all $X$ in $\real^{m\times m}$
 and $Y$ in $\real^{n\times n}$, we have
 $$
 0=\langle A-E, XE+EY \rangle
 =\langle (A-E)E^T,X \rangle 
 + \langle E^T(A-E), Y \rangle .
 $$
 It follows that $E$ satisfies the 
 given equations.
 
 (iii) is equivalent to (iv): This 
 equivalence is obvious.
 \end{proof}
  
 \begin{prop}{\bf Quasi-commuting relations.}
 Let $E$ be an equilibrium point of the
 differential equation (*). Then
 \begin{itemize}
 \item
 The matrix $E$ satisfies the equations
 $$
 AE^T=EA^T,\quad A^TE=E^TA.
 $$
 \item
 The matrix $E$ satisfies the equations
 $$
 A^TAE=E^TAA^T,\quad AA^TE=EA^TA.
 $$
 \end{itemize}
 \end{prop}
 
 I call the two equations which appear in the 
 first conclusion of this proposition, the 
 ``quasi-commuting'' relations for $E$.
 
 \begin{proof}
 We have $AE^T=EE^T$ and $E^TA=E^TE$ from
 the proposition characterizing the
 equilibrium points. To get the
 quasi-commuting relations we simply
 use the symmetry of $EE^T$ and $E^TE$.
 
 To get the other relations we simply
 apply the quasi-commuting relations repeatedly.
 In particular, we have
 \begin{itemize}
 \item
 $A^T(AE^T)=A^T(EA^T)=(A^TE)A^T=(E^TA)A^T$ and
 \item
 $A(A^TE)=A(E^TA)=(AE^T)A=(EA^T)A$.
 \end {itemize}
 \end{proof}
 
 In the following proof and example, I use
 $E^{pq}$ to denote the $m\times n$ matrix
 with a one in position $pq$ and zeros elsewhere:
 $E^{pq}_{ij}:=\delta(i,p)\delta(j,q)$.
 Note that these matrices form a basis of
 the vector space $\real^{m\times n}$.
 
 \begin{prop}
 {\bf Stability of the equilibrium points.}
 If the matrix $A$ has distinct positive 
 singular values, then the differential 
 equation (*) has isolated equilibrium
 points only one of which is stable. 
 It follows that the solutions of the
 differential equation converge and 
 that almost all of them converge 
 to the stable equilibrium point.
 \end{prop}
 
 \remarkr Note that the set of 
 matrices with distinct positive
 singular values is a generic
 (that is, an open and dense) subset
 of $\real^{m\times n}$.
 
 \begin{proof}
 We do the case $m\geq n$. The proof
 in the case $m\leq n$ is essentially
 the same.
 
 We have been working in a coordinate-free
 way until now. We now choose a 
 convenient coordinate system in 
 which to do calculations. In 
 particular, we choose the basis so that
 $A$ is a diagonal matrix of ordered
 singular values:
 $$
 A=\diag(\sigma_1 > \sigma_2 > \dots >
 \sigma_n).
 $$
 
 \claim If a matrix $E$ is an equilibrium
 point of the differential equation
 then $E$ is a diagonal matrix.
 
 Recall that $E$ must satisfy $AA^TE=EA^TA$. 
 We simply calculate these matrix
 products and compare entries. 
 We have 
 $
 (AA^TE)_{ij}= \sigma_i^2 E_{ij}
 $
 if $i\leq n$ and 
 $(AA^TE)_{ij}=0$
 if $i>n$.
 We also have
 $(EA^TA)_{ij}= E_{ij}\sigma_j^2$.
 We conclude that
 if $i\leq n$ and $i\neq j$ 
 then 
 $\sigma_i^2E_{ij}=E_{ij}\sigma_j^2$
 and hence $E_{ij}=0$ since 
 $\sigma_i^2\neq \sigma_j^2$.
 If $i>n$ and $i\neq j$ then
 $0=E_{ij}\sigma_j^2$ 
 and hence $E_{ij}=0$ since
 $\sigma_j^2\neq 0$.
 Thus all the off-diagonal
 entries of $E$ must be zero. 
 
 \claim Let $E:=\diag(e_1,\dots,e_n)$
 be an equilibrium point of the
 differential equation. Then, for
 $i=1,2,\dots,n$, either 
 $e_i=\sigma_i$ or $e_i=0$.
 
 Since the vector field vanishes at
 $E$, we have 
 $$
 0=(\sigma_i-e_i)e_i^2 + e_i^2(\sigma_i-e_i)
 =2(\sigma_i-e_i)e_i^2.
 $$
 
 \claim The solutions of the differential
 equation converge.
 
 From the last claim we see that there
 are a finite number of equilibrium points. 
 A gradient flow confined to 
 a compact set with a finite number of 
 equilibrium points must converge. See,
 for example, Palis and de Melo(1982).

 We now turn to the classification of 
 the equilibrium points. 
 Let $E=\diag(e_1,\dots,e_n)$ be 
 an equilibrium point.  
 We compute the linearization of the 
 differential equation at $E$: We
 get the linear differential equation 
 $$
 X^\prime=D.F.E.X = (A-E)X^TE+EX^T(A-E)
 - XE^TE - EE^TX.
 $$
 We regard $D.F.E$ as a linear
 map on the space tangent to 
 the orbit of $E$ at $E$.
 (By the way, it is easy to check that
 this map is self-adjoint.)
 The nature of the equilibrium is 
 determined by this linear map. In 
 particular, the equilibrium point $E$
 is stable if the eigenvalues of this
 map are all negative. If this
 map has a positive eigenvalue then
 the equilibrium point is unstable.
 We want to see that exactly one of 
 the equilibrium points has all 
 eigenvalues negative (a stable
 situation) and that all 
 of the other equilibrium points
 have at least one positive eigenvalue 
 (an unstable situation).
   
 We have
 $$
 (D.F.E.X)_{ii}= c_i x_{ii}.
 $$
 where $c_i:=2e_i(\sigma_i-2e_i)$.
 For $i\neq j$, we have
 $$
 (D.F.E.X)_{ij}=
 b_{ij}x_{ji}- a_{ij}x_{ij}
 $$
 where
 $a_{ij}:=(e_i^2 + e_j^2)$
 and
 $b_{ij}:= (\sigma_i-e_i)e_j + e_i(\sigma_j-e_j)$.
 
 Since the $ij$ entry of $D.F.E.X$
 involves only the $ij$ and $ji$ entry
 of $X$, we temporarily restrict our attention
 to 2 by 2 matrices. 
 
 We need to find the eigenvalues
 of the map:
 $$
 \begin{pmatrix}
 x_{ij} \\ x_{ji}
 \end{pmatrix}
 \mapsto
 M 
 \begin{pmatrix}
 x_{ij} \\ x_{ji}
 \end{pmatrix},
 $$
 where 
 $$
 M:=
 \begin{pmatrix}
 - a_{ij} & b_{ij}\\
 b_{ij} & -a_{ij}
 \end{pmatrix}.
 $$
 The matrix $M$ has the 
 following form:
 $$
 \begin{pmatrix}
 -a & b \\
 b  & -a
 \end{pmatrix}.
 $$
  This matrix has eigenvalues
 $-a\pm b$. 
 In particular, 
 $$
 \begin{pmatrix}
 -a & b \\
  b & -a 
 \end{pmatrix}
 \begin{pmatrix}
 1 \\ 1
 \end{pmatrix}
 = (-a+b)
 \begin{pmatrix}
 1 \\ 1
 \end{pmatrix}
 $$
 and
 $$
 \begin{pmatrix}
 -a & b \\
  b & -a 
 \end{pmatrix}
 \begin{pmatrix}
 1 \\ -1
 \end{pmatrix}
 = -(a+b)
 \begin{pmatrix}
 1 \\ -1
 \end{pmatrix}.
 $$
 Hence the eigenvalues of $M$
 are
 $$
 \lambda_1:=
 -(e_i+e_j)^2+\sigma_i e_j + e_i\sigma_j
 $$
 and
 $$
 \lambda_2:=
 -(e_i-e_j)^2-\sigma_i e_j - e_i\sigma_j.
 $$
Note that $\lambda_2\leq 0$ for all values
of $e_i,e_j,\sigma_i$ and $\sigma_j$
since these values are always nonnegative.

\claim The diagonal matrix 
 $E^*:=\diag(\sigma_1,\dots,\sigma_k,0,\dots,0)$,
 where $k$ is the rank of the initial
 matrix $K$, is a stable equilibrium
 point. 

We want to see that all the eigenvalues
associated with this equilibrium point
are negative.
Note that the set 
$\{E^{pq}:1\leq p\leq k \text{ or } 1\leq q \leq k\}$
is a basis of the space $\Tan.\Orbit(E^*).E$
tangent to the orbit of $E^*$ at $E^*$.  

If $1\leq i\leq k$ and $1\leq j\leq k$ and 
$i\neq j$ then
$a_{ij}=-(\sigma_i^2+\sigma_j^2)$;
if $1\leq j\leq k$ and $k<j\leq n$ then
$a_{ij}=-\sigma_i^2$;
if $k<i\leq m$ and $1\leq j\leq n$ then
$a_{ij}=-\sigma_j^2$.
If $(1\leq i\leq k$ or $1\leq j\leq k)$
and $i\neq j$ then
$b_{ij}=0$.
For $i=1,\dots,k$,
$c_i=-\sigma_i^2$.
The eigenvalue-vector pairs of $D.F.E^*$
are
\begin{itemize}
\item
$(-(\sigma_i^2+\sigma_j^2),E^{ij})$
for $1\leq i \leq k, 1\leq j\leq k, i\neq j$,
\item
$(-\sigma_i^2,E^{ij})$ 
for $1\leq i\leq k$ and $k<j\leq n$,
\item
$(-\sigma_j^2,E^{ij})$ 
for $k< i\leq m$ and $1 \leq j\leq k$,
\item
$(-\sigma_k^2,E_{ii})$ for $1\leq i \leq k$.
\end{itemize}
Note all these eigenvalues are negative. 

\smallskip
\claim If $E$ is an equilibrium
 point is different than $E^*$
 then $E$ is unstable.
 
In this case the set
$\{E^{pq}: e_p\neq 0 \text{ or } e_q\neq 0\}$
is a basis for the tangent space. 

We want to see that the linear
map $D.F.E$ on the tangent space has
at least one positive eigenvalue. 
Since $E$ is different than $E^*$,
there is an index $p$ satisfying
$1\leq p \leq k$ and $e_p=0$.
Since $E$ has rank $k$, there
is an index $q$ satisfying 
$p<q$ and $e_q\neq 0$. Then
$e_q=\sigma_q$. Note that $E^{pq}$
and $E^{qp}$ are in the tangent space. 
We have
$$
D.F.E.(E^{pq}+E^{qp})
= (b_{pq}+ a_{pq})(E^{pq}+E^{qp})
= (\sigma_p-\sigma_q)\sigma_q(E^{pq}+E^{qp}).
$$
 \end{proof}
  
 \example We do the case $m=4, n=3$ 
 to illustrate the calculations which
 appear in the proof of the last proposition. 
 We consider
 $A:=\diag(\sigma_1,\sigma_2,\sigma_3)$
 where 
 $\sigma_1> \sigma_2 > \sigma_3> 0$.
 Let 
 $$
 E:=
 \begin{pmatrix}
 e_{11} & e_{12} & e_{13}\\
 e_{21} & e_{22} & e_{23}\\
 e_{31} & e_{32} & e_{33}\\
 e_{41} & e_{42} & e_{43}
 \end{pmatrix}
$$
be an equilibrium point.
We have that $AA^T$ is the 4 by 4 diagonal matrix 
$\diag(\sigma_1^2,\sigma_2^2,\sigma_3^2,0)$
and $A^TA$ is the 3 by 3 diagonal matrix
$\diag(\sigma_1^2,\sigma_2^2,\sigma_3^3)$. Hence
$$
AA^TE =
\begin{pmatrix}
 \sigma_1^2 e_{11} & \sigma_1^2 e_{12} 
 &\sigma_1^2 e_{13}\\
\sigma_2^2 e_{21} & \sigma_2^2 e_{22} 
 &\sigma_2^2 e_{13}\\
\sigma_3^2 e_{31} & \sigma_3^2 e_{32}
 &\sigma_3^2 e_{33}\\
 0 & 0 &0
 \end{pmatrix}
$$
and 
$$
EA^TA =
 \begin{pmatrix}
\sigma_1^2 e_{11} & \sigma_2^2 e_{12} 
 &\sigma_3^2 e_{13}\\
\sigma_1^2 e_{21} & \sigma_2^2 e_{22} 
 &\sigma_3^2 e_{23}\\
\sigma_1^2 e_{31} & \sigma_2^2 e_{32} 
 &\sigma_3^2 e_{33}\\
\sigma_1^2 e_{41} & \sigma_2^2 e_{42} 
 &\sigma_3^2 e_{43}\\
 \end{pmatrix}.
 $$
 Equating the entries of these two
 matrices, we see that 
 that all the off-diagonal
 entries of $E$ must be zero. 
 
 We now set $E:=\diag(e_1,e_2,e_3)$.
 We consider the equilibrium equation
 $AE^T=EE^T$. We have
 $AE^T=\diag(\sigma_1 e_1, \sigma_2 e_2, \sigma_3 e_3,0)$
 and 
 $EE^T=\diag(e_1^2,e_2^2,e_3^2,0)$. Equating
 the entries of these two matrices we 
 get, for i=1,2,3,
 $
 \sigma_i e_i = e_i^2,
 $
and hence 
 $e_i=\sigma_i$ or $e_i=0$.
 The specified low rank $k$ will determine the
 number of $e_i$ which are zero.
 
 We turn to the stability classification of the
 equilibrium points. We have
 \begin{align*}
 D.F.E.X &=(A-E)X^TE + EX^T(A-E) -XE^TE - EE^TX \\ 
&=
 \begin{pmatrix}
   (\sigma_1-e_1)x_{11}e_1 
 & (\sigma_1-e_1)x_{21}e_2
 & (\sigma_1-e_1)x_{31}e_3 \\
   (\sigma_2-e_2)x_{12}e_1 
 & (\sigma_2-e_2)x_{22}e_2
 & (\sigma_2-e_2)x_{32}e_3 \\
   (\sigma_3-e_3)x_{13}e_1 
 & (\sigma_3-e_3)x_{23}e_2
 & (\sigma_3-e_3)x_{33}e_3 \\
 0 &0 &0
 \end{pmatrix} \\
 &+
 \begin{pmatrix}
   e_1x_{11}(\sigma_1-e_1) 
 & e_1x_{21}(\sigma_2-e_2)
 & e_1x_{31}(\sigma_3-e_3) \\
   e_2x_{12}(\sigma_1-e_1)
 & e_2x_{22}(\sigma_2-e_2)
 & e_2x_{32}(\sigma_3-e_3) \\
   e_3x_{13}(\sigma_1-e_1) 
 & e_3x_{23}(\sigma_2-e_2)
 & e_3x_{33}(\sigma_3-e_3) \\
 0 &0 &0
\end{pmatrix} \\
&-\begin{pmatrix}
   x_{11}e_1^2 & x_{12}e_2^2 & x_{13}e_3^2 \\
   x_{21}e_1^2 & x_{22}e_2^2 & x_{23}e_3^2 \\
   x_{31}e_1^2 & x_{32}e_2^2 & x_{33}e_3^2 \\
   x_{41}e_1^2 & x_{42}e_2^2 & x_{43}e_3^2
   \end{pmatrix}
-\begin{pmatrix}
   e_1^2x_{11} & e_1^2x_{21} & e_1^2x_{31} \\
   e_2^2x_{12} & e_2^2x_{22} & e_2^2x_{32} \\
   e_3^2x_{13} & e_3^2x_{23} & e_3^2x_{33} \\
   0 &0 &0
   \end{pmatrix}
   \\
&=
\begin{pmatrix}
   0           & b_{12}x_{21}& b_{13}x_{31} \\
   b_{21}x_{12}& 0           & b_{23}x_{32} \\
   b_{31}x_{13}& b_{32}x_{23}& 0            \\
 0 &0 &0
 \end{pmatrix}
 +
 \begin{pmatrix}
 c_1 x_{11}    & -a_{12}x_{12} & -a_{13}x_{13} \\
 -a_{21}x_{21} & c_2 x_{22}    & -a_{23}x_{23} \\
 -a_{31}x_{31} & -a_{32}x_{32} & c_3 x_{33}\\
 -e_1^2x_{41}  & -e_2^2x_{42}  & -e_3^2x_{43}
 \end{pmatrix}
 \end{align*}
 where (using the same notation as that of the 
 proof)
 \begin{align*}
 a_{ij}&:= e_i^2 + e_j^2, \\
 b_{ij}&:= (\sigma_i-e_i)e_j + e_i(\sigma_j-e_j), \\
 c_i &:= 2e_i(\sigma_i - 2e_i).
 \end{align*}

We now consider the stability of the
equilibrium points when $k:=2$ is the
given rank. There are three cases.

\Case $E:=E^*:=\diag(\sigma_1,\sigma_2,0)$

Note that the tangent space to 
$\Orbit(E)$ at $E$ consists of matrices 
$VE+EW$ where $V$ is in $\real^{4\times 4}$
and $W$ is in $\real^{3\times 3}$. It is
easy to see that these matrices have 
the following form:
$$
X:=
\begin{pmatrix}
x_{11} & x_{12} & x_{13} \\
x_{21} & x_{22} & x_{23} \\
x_{31} & x_{32} & 0 \\
x_{41} & x_{42} & 0 
\end{pmatrix}.
$$
If $X$ is such a matrix then
$$
D.F.E^*.X =
\begin{pmatrix}
-2\sigma_1^2 x_{11} 
& -(\sigma_1^2+\sigma_2^2)x_{12} 
& -\sigma_1^2 x_{13} \\
-(\sigma_1^2+\sigma_2^2)x_{21} 
& -2\sigma_2^2 x_{22} 
& -\sigma_2^2 x_{23} \\
-\sigma_1^2 x_{31} 
& -\sigma_2^2 x_{32} 
& 0 \\
-\sigma_1^2 x_{41} 
& -\sigma_2^2 x_{42} 
& 0 
\end{pmatrix}
$$
since
\begin{align*}
a_{12}&= \sigma_1^2 + \sigma_2^2,
a_{13}= \sigma_1^2,
a_{23}= \sigma_2^2, \\ 
b_{12}&= b_{13}= b_{23}=0, \\
c_1 &= -2\sigma_1^2,
c_2 = -2\sigma_2^2, \text{ and } 
c_3 = 0 .
\end{align*}
The eigenvalues of $D.F.E^*$ are
all strictly negative.
In particular, the eigenvalue-vector
pairs are
\begin{align*}
&(-2\sigma_1^2,E^{11}),
(-(\sigma_1^2+\sigma_2^2),E^{12}),
(-\sigma_1^2,E^{13}), \\
&(-(\sigma_1^2+\sigma_2^2),E^{21}),
(-2\sigma_2^2,E^{22}),
(-\sigma_2^2,E^{23}), \\
&(-\sigma_1^2,E^{31}),
(-\sigma_2^2,E^{32}), \\
&(-\sigma_1^2,E^{41}),
(-\sigma_2^2,E^{42}).
\end{align*}
 \Case $E:=\diag(\sigma_1,0,\sigma_3)$

Then the tangent space to $\Orbit(E)$ at $E$
consists of matrices having the following
form:
$$
X:=
\begin{pmatrix}
x_{11} & x_{12} & x_{13} \\
x_{21} & 0      & x_{23} \\
x_{31} & x_{32} & x_{33} \\
x_{41} & 0      & x_{43}
\end{pmatrix}.
$$
If $X$ is such a matrix then $D.F.E.X$
is the following matrix:
$$
\begin{pmatrix}
0
& \sigma_1 \sigma_2 x_{21} 
&0 \\
\sigma_1 \sigma_2 x_{12} 
&0 
& \sigma_2 \sigma_3 x_{32} \\
0 
&\sigma_2 \sigma_3 x_{23} 
&0 \\
0 & 0 & 0
\end{pmatrix}
+
\begin{pmatrix}
-2\sigma_1^2 x_{11} 
& -\sigma_1^2 x_{12} 
& -(\sigma_1^2+\sigma_3^2) x_{13} \\
-\sigma_1^2 x_{21} 
& 0
& -\sigma_3^2 x_{23} \\
-(\sigma_1^2+\sigma_3^2) x_{31} 
& -\sigma_3^2 x_{32} 
& -2\sigma_3^2 x_{33} \\
-\sigma_1^2 x_{41} 
& 0
& -\sigma_3^2 x_{43} 
\end{pmatrix}
$$
since
\begin{align*}
a_{12}&= \sigma_1^2,
a_{13}= \sigma_1^2 + \sigma_3^2,
a_{23}= \sigma_3^2, \\ 
b_{12}&= \sigma_1\sigma_2,
b_{13}= 0,
b_{23}=\sigma_2\sigma_3, \\
c_1 &= -2\sigma_1^2,
c_2 = 0, \text{ and } 
c_3 = -2\sigma_3^2.
\end{align*}
There is a positive eigenvalue. In 
particular, 
$$
D.F.E.(E^{23}+E^{32})
= (\sigma_2-\sigma_3)\sigma_3(E^{23}+E^{32}).
$$
 
\Case $E:=\diag(0,\sigma_2,\sigma_3)$

Then the tangent space to $\Orbit(E)$ at $E$
consists of matrices having the following
form:
$$
X:=
\begin{pmatrix}
0      & x_{12} & x_{13} \\
x_{21} & x_{22} & x_{23} \\
x_{31} & x_{32} & x_{33} \\
0      & x_{42} & x_{43}
\end{pmatrix}.
$$
If $X$ is such a matrix then
$D.F.E.X$ is 
$$
\begin{pmatrix}
0 
& \sigma_1 \sigma_2 x_{21} 
& \sigma_2\sigma_3 x_{31} \\
\sigma_1 \sigma_2 x_{12}
&0 
&0 \\
\sigma_1 \sigma_3 x_{13} 
&0
&0 \\
0 & 0 & 0
\end{pmatrix}
+
\begin{pmatrix}
0
& -\sigma_2^2 x_{12} 
& -\sigma_3^2 x_{13} \\
-\sigma_2^2 x_{21} 
& -2\sigma_2^2 x_{22} 
& -(\sigma_2^2+\sigma_3^2) x_{23} \\
-\sigma_3^2 x_{31} 
& -(\sigma_2^2+\sigma_3^2) x_{32} 
& -2\sigma_3^2 x_{33} \\
0
& -\sigma_2^2 x_{42} 
& -\sigma_3^2 x_{43} 
\end{pmatrix}
$$
since
\begin{align*}
a_{12}&= \sigma_2^2,
a_{13}= \sigma_3^2,
a_{23}= \sigma_2^2+\sigma_3^2, \\ 
b_{12}&= \sigma_1 \sigma_2,
b_{13}= \sigma_1 \sigma_3,
b_{23}=0, \\
c_1 &=0,
c_2 = -2\sigma_2^2, \text{ and } 
c_3 = -2 \sigma_3^2.
\end{align*}
There is a positive eigenvalue.
In particular, 
$$
D.F.E.(E^{12}+E^{21})
= (\sigma_1-\sigma_2)\sigma_2(E^{12}+E^{21}).
$$

\newpage
\section{Acknowledgements}

I did this research during the fall
of 2006 while visiting the Institute
for Mathematics and its Applications
(IMA) at the University of 
Minnesota. I thank the members of 
the IMA for their hospitality
and support. They have created
a very stimulating intellectual
environment.  
Greg Reid (Department of 
Applied Mathematics, University of
Western Ontario) and Wenyuan Wu
(Department of Applied Mathematics,
University of Western Ontario),
who were at the IMA during the
fall of 2006, suggested that I 
look at this problem; I thank them 
for their encouragement. I thank Dong E.
Chang (Department of Applied Mathematics,
University of Waterloo) and  
Philip Rostalski (Automatic
Control Laboratory, ETH Zurich)
for their careful reading and valuable comments
of an early draft of this paper.   
I thank Henry Wolkowitz (Department
of Combinatorics and Optimization, 
University of Waterloo) for his 
advice and encouragement. 

\newpage

\section{Appendix: The Frobenius Inner Product }

We use the ``Frobenius'' (or ``euclidean'')
inner product in the space $\real^{m\times n}$
of $m$ by $n$ real matrices.
For $X$ and $Y$ in this space, the 
{\bf Frobenius inner product} is 
defined by
$$
\langle X,Y \rangle := \Trace(XY^T).
$$
In terms of coordinates, 
$\langle X,Y \rangle =
\sum\{X_{ij}Y_{ij}:i=1,\dots,m,
j=1,\dots,n\}$.
Here we review a few of the properties of 
this inner product.

\begin{prop}{\bf Adjoints of 
multiplication maps.}
Let $B$ and $Z$ be elements 
of $\real^{m\times n}$.
\begin{itemize}
\item
For $X\in\real^{m\times m}$,
$
\langle XB,Z \rangle = 
\langle X,ZB^T \rangle.
$
\item
For $Y\in\real^{n\times n}$,
$
\langle BY,Z \rangle = 
\langle Y,B^TZ \rangle.
$
\end{itemize}
\end{prop}

\begin{proof}
We have
$$
\Trace(XBZ^T)=\Trace(X(ZB^T)^T)
$$
and
$$
\Trace(BYZ^T)=\Trace(YZ^TB)=
\Trace(Y(B^TZ)^T).
$$
\end{proof}

\begin{prop}
{\bf Orthogonal invariance.}
Let $U$ be an $m$ by $m$ real orthogonal matrix
and let $V$ be an $n$ by $n$ real orthogonal 
matrix. Then, for all $X$ and $Y$ 
in $\real^{m\times n}$, 
\begin{itemize}
\item
$\langle UX,UY \rangle = \langle X,Y \rangle$ and
\item
$\langle XV,YV \rangle = \langle X,Y \rangle$.
\end{itemize}
\end{prop}

\begin{proof}
By the result concerning the adjoints of
multiplication maps, we have:
\begin{itemize}
\item
$\langle UX,UY \rangle
= \langle X, U^TUY \rangle
= \langle X,Y \rangle$ and
\item
$\langle XV,YV \rangle
= \langle X,YVV^T \rangle
= \langle X,Y \rangle$
\end{itemize}
\end{proof}

\newpage
\section{References}

\begin{itemize}
\item
Birkhoff, G. and MacLane, S. (1953)
{\it A Survey of Modern Algebra}, 
Macmillan.
\smallskip

\item
Chu, M. (1984)
{\it The generalized Toda flow, the
QR algorithm, and the centre manifold
theory},
SIAM J. Alg. Discr. Math. 5, 187-201.
\smallskip

\item
Chu, M. (1986a)
{\it A differential equation approach to
the singular value decomposition
of bidiagonal matrices},
Lin. Alg. Appl. 80, 71-80.
\smallskip

\item
Chu, M. (1986b)
{\it A continuous approximation to the
generalized Schur decomposition},
Lin. Alg. Appl. 78, 119-132.
\smallskip

\item
Chu, M. and Driessel, K.R. (1990)
{\it The projected gradient method
for least squares approximation with
spectral constraints},
SIAM J. Numerical Analysis 27, 1050-1060.
\smallskip

\item
Deift, P., Nanda, T. and Tomei, C. (1983)
{\it Differential equations for the 
symmetric eigenvalue problem},
SIAM J. Numer. Analysis 20, 1-22.
\smallskip

\item
Demmel, J.W. (1997)
{\it Applied Numerical Linear Algebra},
SIAM.
\smallskip

\item
Driessel, K.R. (2004) 
{\it On computing cannonical forms using flows}, 
Lin. Alg. Appl. 379, 353-379.
\smallskip

\item
Driessel, K.R. and Gerisch, A.(2007) 
{\it Zero-preserving iso-spectral
flows bases on parallel sums}, 
Lin. Alg. Appl. 421, 69-84.
\smallskip

\item
Eckart, G. and Young, G.(1936) 
{\it The approximation of one
matrix by another of lower rank}, 
Psychometrika 1, 221-218.
\smallskip


\item
Helmke, U. and Moore, J.B. (1995)
{\it Optimization and Dynamical Systems},
Springer.
\smallskip

\item
Helmke, U. and Shayman, M.A. (1995)
{\it Critical points of matrix least squares
distance functions},
Lin. Alg. Appl. 215, 1-19.
\smallskip

\item
Hirsch, M.W. and Smale, S. (1974)
{\it Differential Equations, Dynamical Systems,
and Linear Algebra},
Academic Press.
\smallskip

\item
Horn, R.A. and Johnson, C.R. (1985)
{\it Matrix Analysis},
Cambridge University Press.
\smallskip


\item
Nanda, T. (1982) 
{\it Isospectral flows on band matrices},
Doctoral Dissertation, Courant Institute, 
New York. 
\smallskip

\item
Nanda, T. (1985) 
{\it Differential equations and the 
QR algorithm},
SIAM J. Numer. Analysis 22, 310-321.
\smallskip

\item
Palis, J., Jr. and de Melo, W. (1982) 
{\it Geometric Theory of 
Dynamical Systems},
Springer.
\smallskip

\item 
Stewart, G.W. (1993)
{\it On the early history of the 
singular value decomposition},
SIAM Review 35, 551-566.
\smallskip

\item
Symes, W.W. (1980a)
{\it Systems of Toda type, inverse 
spectral problems, and representation
theory}, 
Inventiones Mathematicae 59, 13-51.
\smallskip

\item
Symes, W.W. (1980b)
{\it Hamiltonian group actions and integrable
systems}, 
Physica 1D, 339-374.
\smallskip

\item
Symes, W.W. (1982)
{\it The QR algorithm and scattering for
the finite nonperiodic Toda lattice},
Physica 4D, 275-280.
\smallskip

\item
Thorpe, J.A. (1979)
{\it Elementary Topics in Differential
Geometry},
Springer.
\smallskip

\item
Watkins, D.S. (1984a) 
{\it Isospectral flows},
SIAM Review 26, 379-392.
\smallskip

\item
Watkins, D.S. (1984b) 
{\it The Toda flow and other isospectral
flows},
Lin. Alg. Appl. 59, 196-201.
\smallskip

\end{itemize}

\end{document}